\documentclass[11pt,twoside,a4paper]{article}
\usepackage{amsmath,amssymb,mathrsfs}
\usepackage{amsthm}
\usepackage[utf8x]{inputenc}
\usepackage{bm}
\usepackage{avant}
\usepackage[english]{babel}
\usepackage{thmtools,thm-restate}
\usepackage[nottoc]{tocbibind}
\usepackage{hyperref}
\usepackage{graphicx}
\usepackage{enumitem}
\usepackage{tikz}
\usepackage{bbm}
\usepackage[inner=2.4cm,outer=2.4cm,top=2.8cm,bottom=2.8cm]{geometry}
\usepackage{mathtools}
\usepackage{makeidx}
\usepackage{bold-extra}
\usepackage{todonotes}

\DeclareMathOperator{\id}{id}

\DeclareMathOperator{\supp}{supp}
\DeclareMathOperator{\dist}{dist}

\newcommand{\scal}[2]{\ensuremath{\langle #1 , #2 \rangle}} 
\newcommand{\norm}[1]{\left\lVert#1\right\rVert}

\newcommand{\Leb}{\mathscr{L}}
\newcommand{\R}{\mathbb{R}}
\newcommand{\Q}{\mathbb{Q}}
\newcommand{\p}{\mathtt p} 
\newcommand{\de}{\ensuremath{\, \mathrm d}} 

\newcommand{\suchthat}{\ensuremath{\,:\,}} 
\newcommand\restr[2]{{
  \left.\kern-\nulldelimiterspace 
  #1 
  \right|_{#2} 
  }}

\DeclareMathOperator{\OptPlans}{OptPlans}
\newcommand{\Prob}{\mathscr{P}}

\usepackage{titlesec}
\titleformat{\section}
  {\scshape\large\centering}
  {\thesection}{0.5em}{}
\titleformat{\subsection}
  {\scshape\centering}
  {\thesubsection}{0.4em}{}

\author{Guoxi Liu\footnote{Trinity College, University of Oxford. \textit{E-mail}:  \href{guoxi.liu@trinity.ox.ac.uk}{guoxi.liu@trinity.ox.ac.uk}}, \ Mattia Magnabosco\footnote{Mathematical Institute, University of Oxford. \textit{E-mail}:  \href{mailto:mattia.magnabosco@maths.ox.ac.uk}{mattia.magnabosco@maths.ox.ac.uk}}, \ Yicheng Xia\footnote{St Anne's College, University of Oxford. \textit{E-mail}:  \href{mailto:yicheng.xia@st-annes.ox.ac.uk}{yicheng.xia@st-annes.ox.ac.uk}}}
\title{\textbf{On the existence of $L^p$-Optimal Transport maps for norms on $\R^N$}}
\date{}

\newtheoremstyle{remark}
        {10pt}
        {10pt}
        {}
        {}
        {\itshape}
        {.}
        {.4em}
        {}

\newtheoremstyle{proof}
        {10pt}
        {10pt}
        {}
        {}
        {\itshape}
        {.}
        {.4em}
        {}
        
\newtheoremstyle{definition}
        {10pt}
        {10pt}
        {}
        {}
        {\bfseries}
        {.}
        {.4em}
        {}

\newtheoremstyle{theorem}
        {10pt}
        {10pt}
        {\slshape}
        {}
        {\bfseries}
        {.}
        {.4em}
        {}

\theoremstyle{theorem}
\newtheorem{theorem}{Theorem}[section]
\newtheorem{prop}[theorem]{Proposition}
\newtheorem{corollary}[theorem]{Corollary}

\theoremstyle{definition}
\newtheorem{definition}[theorem]{Definition}

\theoremstyle{remark}
\newtheorem{remark}[theorem]{Remark}
\newtheorem{example} [theorem]{Example}

\theoremstyle{proof}
\newtheorem*{pro}{Proof}
\newenvironment{pr}{\begin{pro}%
 \pushQED{\qed}}%
 {\popQED\end{pro}}

\makeindex
\begin{document}
\maketitle

\begin{abstract}
    In this paper, we prove existence of $L^p$-optimal transport maps with $p\in (1,\infty)$ in a class of branching metric spaces defined on $\R^N$. In particular, we introduce the notion of \textit{cylinder-like} convex function and we prove an existence result for the Monge problem with cost functions of the type $c(x, y) = f(g(y - x))$, where $f: [0, \infty) \rightarrow [0, \infty)$ is an increasing strictly convex function and $g: \R^N \rightarrow [0, \infty)$ is a cylinder-like convex function. When specialised to cylinder-like norm, our results shows existence of $L^p$-optimal transport maps for several ``branching'' norms, including all norms in $\R^2$ and all crystalline norms.
\end{abstract}

\section{Introduction}

The theory of optimal transport has been deeply studied in recent years due to its many applications in pure and applied mathematics. The optimal transport problem was originally introduced by Monge in the 18th century and its formulation in $\R^N$ is as follows: given two Borel probability measures $\mu,\nu\in\Prob(\R^N)$ (called marginals) and a non-negative Borel cost function $c: \R^N \times \R^N   \to  [0,\infty]$, study the following minimisation problem:
\begin{equation} \label{MongeFormulation}
	\inf\left\{ \int_X c(x,T(x)) \de \mu(x) \suchthat \text{$T:\R^N \to \R^N$ Borel, $T_\# \mu =\nu$} \right\}. \tag{M}
\end{equation}
For several applications it is particularly interesting to identify eventual minimisers of \eqref{MongeFormulation}, which are called optimal transport maps. However, the infimum of the optimal transport problem is not realised in general and therefore optimal transport maps do not exist for many marginals and cost functions. For the development of a general theory it was thus necessary a breakthrough due to Kantorovich in 1942. He introduce a more general formulation of the optimal transport problem for which the minimisers exist in many cases. In particular, defining the set of admissible transport plans from $\mu$ to $\nu$ as
\begin{equation*}
\mathsf{Adm}(\mu,\nu) := \{ \pi\in\Prob(\R^N\times \R^N) \suchthat (\p_1)_\#\pi = \mu \,\,\text{and} \,\, (\p_2)_\#\pi = \nu\},
\end{equation*}
the Kantorovich's formulation of the optimal transport problem aims to study the minimisation:
\begin{equation} \label{KantorovichFormulation}
	\inf \left\{ \int_{X\times Y} c(x,y) \de \pi(x,y) \suchthat \pi\in\mathsf{Adm}(\mu,\nu) \right\}. \tag{K}
\end{equation}
In this version of the optimal transport problem the minimum is attained whenever the cost function $c:\R^N\times \R^N \to [0,\infty]$ is lower semicontinuous. The minimisers of the Kantorovich formulation are called optimal transport plans and the set containing them all will be denoted by $\OptPlans(\mu,\nu)$. An admissible transport plan $\pi\in\mathsf{Adm}(\mu,\nu) $ is said to be induced by a map if there exists a $\mu$-measurable map $T:\R^N \to \R^N$ such that $\pi=(\id,T)_\# \mu$. Notice that these admissible plans are the ones considered in the Monge Formulation \eqref{MongeFormulation} and thus if $\pi=(\id,T)_\# \mu$ is an optimal transport plan then $T$ is an optimal transport map. 

The problem of addressing existence and uniqueness of optimal transport maps has been already deeply studied. The first result in this direction was obtained by Brenier in \cite{brenier} for the $L^2$-optimal transport problem in $\R^N$, i.e. with cost function $c(x,y)=\norm{x-y}^2_{eu}$. He proved that, whenever the first marginal $\mu$ is absolutely continuous with respect to the Lebesgue measure $\Leb^N$, there exists a unique optimal transport plan in $\OptPlans(\mu,\nu)$ and it is induced by a map. Later, Brenier's result was generalised for the $L^2$-optimal transport problem, i.e. with cost $c(x,y)=d(x,y)^2$, in smooth and nonsmooth metric measure spaces $(X,d,\mathfrak m)$. In particular, for smooth spaces, we mention the works by McCann \cite{McCann} in the Riemannian setting and by Ambrosio-Rigot \cite{AmbrosioRigot} and Figalli-Rifford \cite{FigalliRifford} in the sub-Riemannian setting. In nonsmooth spaces, existence and uniqueness of $L^2$-optimal transport maps is usually proved under a synthetic curvature bound and a non-branching assumption. We mention for example the following articles: \cite{Bertrand,RajalaSchultz} for Alexandrov spaces, \cite{Gigli, RajalaSturm,GRS,cavallettimondino,Schultz,MagnaboscoRigoni} under various curvature-dimension condition and \cite{CavallettiHuesmann,Kell} under quantitative properties on the reference measure. We remark that, in the nonsmooth setting, it is particularly interesting to study the $L^2$-optimal transport problem, as it is the fundamental tool to formulate the groundbreaking theory of $\mathsf{CD}(K,N)$ spaces, pioneered by Sturm \cite{sturmI, sturmII} and Lott-Villani \cite{lott-villani}.

For the purpose of this work it is fundamental to highlight that, as shown by Rajala in \cite{Rajala} (see also \cite{Magnabosco1,Magnabosco2}), a non-branching assumption is necessary to prove uniqueness of the optimal transport map. Indeed, even in the simplest examples of branching spaces, such as $(\R^N, \norm{\cdot}_\infty, \Leb^N)$, we cannot hope to have a unique optimal transport map between two given marginals. In this paper, we aim at proving existence of $L^p$-optimal transport maps with $p\in (1,\infty)$, i.e. for the cost $c(x,y)=d(x,y)^p$, in a class of branching spaces defined on $\R^N$. More specifically, we introduce the notion of \textit{cylinder-like} convex function and we prove an existence result (Theorem \ref{theorem: CylinderLike}) for cost functions of the type $c(x, y) = f(g(y - x))$, where $f: [0, \infty) \rightarrow [0, \infty)$ is an increasing strictly convex function and $g: \R^N \rightarrow [0, \infty)$ is a cylinder-like convex function. Roughly speaking, a convex function is cylinder-like if its level sets look like the curved surface of a cylinder. In particular, we obtain existence of optimal transport maps for the $L^p$ problem in $\R^N$ equipped with a cylinder-like norm. This class of norms contains a good amount of ``branching'' norms, including all norms in $\R^2$ and all crystalline norms. We remark that, for the latter, existence of $L^2$-optimal transport maps was firstly proved in \cite{magnabosco3}.

The strategy we use to prove our main result consists in considering a secondary variational optimisation problem, obtained by minimising a second energy functional among all optimal transport plans. Then, proving that all minimisers of the secondary variational problem are induced by a map, will be sufficient to prove existence of optimal transport maps. This strategy was firstly applied in branching spaces by Rajala \cite{Rajala} (see also \cite{magnabosco3,LMT}), but it has been successful even in proving existence of $L^1$-optimal transport maps, see for example \cite{AmbKirPra, ChaDeP, Pingetal}. We highlight that the $L^1$-optimal transport problem is, for some aspects, similar to the $L^p$-problem in branching spaces, as in both cases the cost function fails to be strictly convex. However, the latter seems more complicated as it is more difficult to identify and deal with the set where the cost function is not strictly convex.

A fundamental notion for our strategy is \emph{cyclical monotonity}, which is one of the most basic and important tools in the theory of optimal transport. We recall that, given a cost function $c: \R^N \times \R^N   \to  [0,\infty]$, a set $\Gamma\subset \R^N\times \R^N$ is said to be $c$-cyclically monotone if 
\begin{equation*}
    \sum_{i=1}^{n} c\left(x_{i}, y_{\sigma(i)}\right) \geq \sum_{i=1}^{n} c\left(x_{i}, y_{i}\right)
\end{equation*}
for every $n\geq1$, every permutation $\sigma$ of $\{1,\dots,n\}$ and every $(x_i,y_i)\in \Gamma$ for $i=1,\dots,n$. As shown by the next proposition, $c$-cyclical monotonicity and optimality are strictly related.

\begin{prop}\label{prop:cmonotonicity}
Let $c: \R^N \times \R^N   \to  [0,\infty]$ be a lower semicontinuous cost function. Then, every optimal transport plan $\pi\in \OptPlans(\mu,\nu)$ such that $\int c \de \pi<\infty$ is concentrated in a $c$-cyclically monotone set. 
\end{prop}

\noindent In Section \ref{sec:secvar}, we recall a ``second-order'' analogous of Proposition \ref{prop:cmonotonicity} for the secondary variational optimisation problem, which will be the central tool in our strategy.

\subsection*{Acknowledgments} 
 G.L. acknowledges support from the Trinity College (University of Oxford) that funded the summer research project. M.M. acknowledges support from the Royal Society through the Newton International Fellowship (award number: NIF$\backslash$R1$\backslash$231659). Y.X. acknowledges support from the Mathematical Institute (University of Oxford) through the departmental funding for summer research projects.

\section{The set of strict convexity}

In this section we introduce the \emph{set of strict convexity} of a given convex function on $\R^N$. This object will play a relevant role in the proofs of our results.
\begin{definition}
 Given a convex function $h : \R^N \to [0, \infty)$, we define its set of strict convexity as 
\begin{equation*}
H_h = \left\{(x,y)\in \R^N \times \R^N: h(tx+(1-t)y)<th(x)+(1-t)h(y),\forall t \in (0,1)\right\}.
\end{equation*}
\end{definition}

\begin{remark}
    By properties of convex functions we see that
    \begin{equation*}
    H_h^c = \left\{(x,y)\in \R^N \times \R^N: h(tx+(1-t)y)=th(x)+(1-t)h(y),\forall t \in (0,1)\right\}.
    \end{equation*}
    Moreover, since any real-valued convex function defined on $\R^N$ is continuous, it's easy to show that $H_h^c$ contains all its limit points, and so is closed. As a consequence, $H_h$ is open.
\end{remark}

\begin{theorem} \label{theorem: StrictlyConvexCase}
Let $h: \R^N \rightarrow [0, \infty)$ be a convex function and consider the cost function $c(x, y) = h(y - x)$. Let $\mu, \nu \in \Prob(\R^N)$ be Borel probability measures on $\R^N$ such that $\mu \ll \Leb^N$ and $\mathscr C(\mu, \nu) < \infty$. Let $\pi \in \OptPlans(\mu, \nu)$ be an optimal plan, and let $(\pi_z)_{z \in \R^N}$ be the disintegration kernel of $\pi$ with respect to the projection $\p_1$ on the first coordinate. For a fixed $z \in \R^N$, define the probability measure $\tilde{\pi}_z$ by $\tilde{\pi}_z(B) = \pi_z(B + z)$ for all $B \subseteq \R^N$ Borel. Then, the set 
\begin{equation*}
E = \left\{ z \in \R^N \colon \exists x, y \in \supp\tilde{\pi}_z \text{ such that } (x, y) \in H_h \right\}
\label{M1}
\end{equation*}
has zero $\mu$-measure.

\end{theorem}

\begin{pr}
Suppose by contradiction that $\mu(E)>0$.
Define the measure $\gamma$ on $\R^N \times \R^N$ as
\begin{equation*}
    \gamma := \int_{\R^N} \tilde{\pi}_z \times \tilde{\pi}_z\, \de \mu(z).
\end{equation*}
 We claim that $\gamma(H_h) > 0$. Indeed, using the fact that $H_h$ is an open subset of $\R^N \times \R^N$ and defining 
\begin{equation*}
E_{a,b,\delta} := \left\{z\in\R^N : \tilde{\pi}_z(B_\delta(a)) > 0, \tilde{\pi}_z(B_\delta(b)) > 0\right\},
\end{equation*}
we have that
\begin{equation*}
    E \subseteq \bigcup_{\substack{a, b \in \mathbb{Q}^N, \,\,  \delta \in \Q^+ \\ B_\delta(a) \times B_\delta(b) \subseteq H_h}} E_{a,b,\delta}.
\end{equation*}
Then, we can find $a, b\in \R^N$ and $\delta>0$ such that $\gamma(H_h) \geq \gamma(B_\delta(a) \times B_\delta(b)) > 0$. Moreover, using the inner regularity of finite Borel measures on $\R^N$, we deduce that there exists a closed set $H_1 \subseteq H_h$ such that $\gamma(H_1)>0$.

Take $(\bar{x},\bar{y}) \in \supp \gamma\big|_{H_1}$. Since $H_1$ is closed, we have that $(\bar{x}, \bar{y}) \in H_1 \subseteq H_h$. Then, consider $\bar{v} = \bar{x} - \bar{y}$, and fix $0 < \epsilon < \frac{1}{2}$. By strict convexity in \(B_\delta(a) \times B_\delta(b)\subseteq H_h\), we have that
\begin{equation*}
    h(\bar{x} - \epsilon \bar{v}) - h(\bar{x}) < h(\bar{y}) - h(\bar{y} + \epsilon \bar{v}).
\end{equation*}
By continuity of $h$, there exists $\delta > 0$ such that
\begin{equation*}
h(x - \epsilon \bar{v}) - h(x) < h(y) - h(y + \epsilon \bar{v})
\quad \text{for all} \quad (x, y) \in B_\delta(\bar{x}) \times B_\delta(\bar{y}).
\end{equation*}
Now, since $(\bar{x}, \bar{y}) \in \supp \gamma\big|_{H_1} \subseteq \supp \gamma$, we deduce that \(\gamma(B_\delta(\bar{x}) \times B_\delta(\bar{y})) > 0\) and consequently  that
\begin{equation*}
F := \left\{ z \in \R^N : \tilde{\pi}_z(B_\delta(\bar{x})) > 0, \; \tilde{\pi}_z(B_\delta(\bar{y})) > 0 \right\}
\end{equation*}
is a set of positive $\mu$-measure. 

On the other hand, by Theorem~\ref{prop:cmonotonicity}, $\pi$ is concentrated on a $c$-cyclically monotone set $\Gamma$. For $z \in \R^N$ write $\Gamma_z := \{y \in \R^N : (z, y) \in \Gamma\}$. Then, we have $\pi_z(\Gamma_z) = 1$ $\mu$-almost everywhere, and thus
\begin{equation*}
G := \left\{ z \in \R^N : \tilde{\pi}_z\big(B_\delta(\bar{x}) \cap (\Gamma_z - z)\big) > 0, \; \tilde{\pi}_z\big(B_\delta(\bar{y}) \cap (\Gamma_z - z)\big) > 0 \right\}
\end{equation*}
is a set of positive $\mu$-measure, hence of positive Lebesgue measure.
By Lebesgue's density point theorem, we can find a density point $\bar{z}$ of $G$. In a neighbourhood of $\bar{z}$ we thus can find $z, z+\eta\bar{v} \in G$, for a suitable $0 < \eta < \epsilon$. Using the definition of $G$, there are \(x\in B_\delta(\bar{x}), y\in B_\delta(\bar{y})\) such that $(z,z+x),(z+\eta\bar{v},z+\eta\bar{v}+y)\in \Gamma$. Therefore, we conclude that
\begin{align*}
     [c(z,z+y+\eta\bar{v}) &+ c(z+\eta\bar{v},z+x)] - [c(z,z+x) + c(z+\eta\bar{v},z+\eta\bar{v}+y)] \\
    & = [h(y+\eta\bar{v}) + h(x - \eta\bar{v})] - [h(x) + h(y)] \\
    & = h\left(\frac{\eta}{\epsilon}(y+\epsilon\bar{v})+\frac{\epsilon-\eta}{\epsilon}(y)\right) + h\left(\frac{\eta}{\epsilon}(x-\epsilon\bar{v})+\frac{\epsilon-\eta}{\epsilon}(x)\right) - h(x) - h(y) \\
    & \leq \left[\frac{\eta}{\epsilon}h(y+\epsilon\bar{v})+\frac{\epsilon-\eta}{\epsilon}h(y)\right] + \left[\frac{\eta}{\epsilon}h(x-\epsilon\bar{v})+\frac{\epsilon-\eta}{\epsilon}h(x)\right] - h(x) - h(y) \\
    & = \frac{\eta}{\epsilon}\left(h(y + \epsilon\bar{v}) + h(x - \epsilon \bar{v}) - h(x) - h(y)\right) < 0.
\end{align*}
This contradicts the $c$-cyclical monotonicity of $\Gamma$.
\end{pr}

As a consequence of the last theorem, we obtain the following corollary, which was already proved by Gangbo and McCann \cite{GangboMcCann} under some additional growth assumption on the cost function. In particular, we provide a simpler proof of an easier result that does not give any structural information about the optimal transport map.

\begin{corollary}
    Let $h: \R^N \rightarrow [0, \infty)$ be a strictly convex function and consider $c(x, y) = h(y - x)$. Let $\mu, \nu \in \Prob(\R^N)$ be Borel probability measures on $\R^N$ such that $\mu \ll \Leb^N$ and $\mathscr C(\mu, \nu) < \infty$.
    Then $\OptPlans(\mu, \nu)$ has a unique element and it is induced by a map.
    \label{cor:strictconv}
\end{corollary}

\begin{pr}
    By assumption there exists at least one optimal transport plan; to show that it is unique and it is induced by a map, it suffices to show that every optimal plan is induced by a map. Indeed, suppose there are two distinct optimal plans $\pi_1, \pi_2 \in \OptPlans (\mu, \nu)$, induced by maps $T_1, T_2$ respectively, which differ on a set of positive $\mu$ measure. Then $(\pi_1 + \pi_2)/2$ is an optimal plan which is not induced by a map, giving a contradiction.

 Hence, take any $\pi \in \OptPlans(\mu, \nu)$. Since $h$ is strictly convex, we have
    \begin{equation*}
        H_h = \{(x, y) \in \R^N \times \R^N : x \neq y\}.
    \end{equation*}
    Thus by Theorem~\ref{theorem: StrictlyConvexCase}, for $\mu$-almost every $z \in \R^N$, $\supp \tilde{\pi}_z$ contains one element. Hence $\pi_z$ is a delta measure for $\mu$-almost all $z \in \R^N$, and so $\pi$ is induced by a map.
\end{pr}

Specialising Corollary \ref{cor:strictconv} to the $L^p$-optimal transport problem for (finite dimensional) normed spaces, we obtain the following result.

\begin{corollary}
    Let $\norm{\cdot}$ be a strictly convex norm on $\R^N$ and $1 < p < \infty$. Consider $c(x, y) = \norm{y - x}^p$ and let $\mu, \nu \in \Prob(\R^N)$ be Borel probability measures with finite $p$th order moment on $\R^N$, such that $\mu \ll \Leb^N$. Then $\OptPlans(\mu, \nu)$ has a unique element and it is induced by a map.
\end{corollary}

\begin{remark}
    Observe that, in the last corollary, the finiteness assumption on the $p$th order moments ensures that $\mathscr C(\mu, \nu) < \infty$.
\end{remark}

\section{Cylinder-like convex functions}

In this section, we introduce the notion of cylinder-like convex function, that will be a central assumption in our main theorem. As already mentioned, the name is due to the fact that, heuristically, a convex function is cylinder-like if its level sets look like the curved surface of a cylinder. In Proposition \ref{CylinderLikeFunctionsProp} and \ref{BuildingCylinderLikeFunctions} we show that the class of cylinder-like functions covers in particular all non-branching and a good amount of branching norms in $\R^N$.

\begin{definition} \label{CylinderLike}

A convex function $g : \R^N \rightarrow [0,\infty)$ is said to be \emph{cylinder-like} if, for every $\bar{x}, \bar{y} \in \R^N$ with $\bar{x} \neq \bar{y}$ and $ g(\bar{x}) = g(\bar{y})$, there exists $\epsilon > 0$ and $\delta > 0$ such that
\begin{equation}\label{eq:CL}
    g(x - \epsilon(\bar{x} - \bar{y})) \leq g(x) \text{   }\,\forall x\in B_\delta(\bar{x})  \quad \text{and} \quad g(y + \epsilon(\bar{x} - \bar{y})) \leq g(y) \text{       }\,\forall y\in B_\delta(\bar{y}),
\end{equation}
where $B_r(z)$ denotes the Euclidean ball of radius $r$ centred in $z$.
\end{definition}

\begin{remark}
    In Definition \ref{CylinderLike}, if \eqref{eq:CL} holds for some $\epsilon > 0$ and $\delta > 0$, then it also holds for every $\varepsilon'\in (0,\varepsilon)$ and $\delta'\in(0,\delta)$.
    \label{CylinderLikeSmallerEpsilonRemark}
\end{remark}

\begin{remark}[Sufficient condition for functions to cylinder-like]
    It suffices to check \eqref{eq:CL} of  Definition~\ref{CylinderLike} for points $\bar{x}, \bar{y} \in \R^N$ with $\bar{x} \neq \bar{y}$, $g(\bar{x}) = g(\bar{y})$ and $ (\bar{x}, \bar{y}) \in H_g^c$.

    To see this, fix $ \epsilon \in (0,1)$ and take any $\bar{x}, \bar{y} \in \R^N$ such that $\bar{x} \neq \bar{y}$, $g(\bar{x}) = g(\bar{y})$ and $(\bar{x}, \bar{y}) \in H_g$. Then, by definition of $H_g$ and since $g(\bar{x}) = g(\bar{y})$, we have
    \begin{equation*}
        g(\bar{x} - \epsilon(\bar{x} - \bar{y})) < g(\bar{x})\quad \text{and}\quad 
        g(\bar{y} + \epsilon(\bar{x} - \bar{y})) < g(\bar{y}).
    \end{equation*}
    Hence, using the continuity of the functions
    \begin{equation*}
        z \mapsto g(z - \epsilon(\bar{x} - \bar{y})) - g(z) \quad \text{and}\quad 
        z \mapsto g(z + \epsilon(\bar{x} - \bar{y}) - g(z),
    \end{equation*}
    we can find a suitable $\delta$ for which \eqref{eq:CL} holds.

\label{CylinderLikeSufficientCondition}
\end{remark}

The next two propositions show that the class of cylinder-like convex functions covers a lot of interesting examples.

\begin{prop}
The following convex functions are cylinder-like:
\begin{itemize}
    \item[(i)] any strictly convex norm on $\R^N$,
    \item[(ii)] any norm on $\R^2$,
    \item[(iii)] the function $f(x) = |\scal{v}{x}|$ where $v \in \R^N$ is any vector.
\end{itemize}

\label{CylinderLikeFunctionsProp}
\end{prop}

\begin{pr}
    \noindent{\slshape (i)} Follows from Remark~\ref{CylinderLikeSufficientCondition}, as given a strictly convex norm $N$ on $\R^N$, no $\bar{x} \neq \bar{y}$ can satisfy $N(x) = N(y)$ and $(\bar{x}, \bar{y}) \in H_N^c$.\vspace{5pt}

    \noindent {\slshape (ii)} Given a norm $\norm{\cdot}$ on $\R^2$, convexity and non-negativity are automatic. By Remark~\ref{CylinderLikeSufficientCondition}, it suffices to check \eqref{eq:CL} (for some $\epsilon > 0$, $\delta > 0$ and) for $\bar{x} \neq \bar{y}$, $\norm{\bar{x}} = \norm{\bar{y}}$ and

    \begin{equation} \label{2DNormFlatPart}
        \norm{t\bar{x} + (1-t)\bar{y}} = \norm{\bar x}= \norm{\bar y} \quad \forall t\in (0,1).
    \end{equation}
    Let $\bar{v} = \bar{x} - \bar{y}$ and define the function
    \begin{equation*}
        f: \R \to [0,+\infty) \qquad f(t) := \norm{\bar{y} + t\bar{v}}.
    \end{equation*} 
    Note that $f$ is convex and constant on $[0, 1]$, by \eqref{2DNormFlatPart}. By the monotonicity of difference quotients, $f$ is increasing on $[1/2, \infty)$ and decreasing on $(-\infty, 1/2]$.
    Define the cones
    \begin{equation*}
    \begin{split}
        &C_1 := \{\lambda(\bar{y} + \mu\bar{v}) : \lambda > 0, \mu \in (1/2, 3/2)\}, \\
        &C_2 := \{\lambda(\bar{y} + \mu\bar{v}) : \lambda > 0, \mu \in (-1/2, 1/2)\}. \\
    \end{split}
    \end{equation*}
    Choose $\delta > 0$ such that
    \begin{equation*}
        \epsilon_0 := \frac{1}{2}\min\left\{\dist(\bar{B}_\delta(\bar{x}), C_1^c), \dist(\bar{B}_\delta(\bar{y}), C_2^c)\right\} > 0,
    \end{equation*}
    where we have used the fact that for $K \subseteq \R^N$ non-empty and compact and $F \subseteq \R^N$ non-empty and closed, $\dist(K, F) := \inf \{d(x, y) : x \in K, y \in F\} > 0$, where $d(x, y)$ denotes the Euclidean distance.

    Take any $x \in B_\delta(\bar{x})$ and call $R:= \norm{\bar v}_{eu}$, where $\norm{\cdot}_{eu}$ denotes the Euclidean norm. Then, we have that $x, x - \frac{\epsilon_0}R \bar{v} \in C_1$. By definition of $C_1$, there exist $\alpha > 0$ and $\beta \in (1/2, 3/2)$ such that
    \begin{equation*}
        x = \alpha (\bar{y} + \beta \bar{v}),
    \end{equation*}
    therefore
    \begin{equation*}
        x - \frac{\epsilon_0}R \bar{v} = \alpha \left(\bar{y} + \left(\beta - \frac{\epsilon_0}{R\alpha}\right)\bar{v}\right).
    \end{equation*}
    As $x, x - \frac{\epsilon_0}R \bar{v} \in C_1$, we have that
    \begin{equation*}
        \frac{1}{2} < \beta - \frac{\epsilon_0}{R\alpha} < \beta < \frac{3}{2}
    \end{equation*}
    and thus, taking into account the monotonicity of $f$ on $[1/2, \infty)$, we conclude that
    \begin{equation}\label{eq:1}
        \norm{x - \frac{\epsilon_0} R \bar{v}} = \alpha f\left(\beta - \frac{\epsilon_0}{R\alpha}\right) \leq \alpha f(\beta) = \norm{x}.
    \end{equation}
    Similarly, using the monotonicity of $f$ on $(-\infty, 1/2]$, we deduce that
    \begin{equation}\label{eq:2}
        \norm{y + \frac{\epsilon_0}{R} \bar{v}} \leq \norm{y}.
    \end{equation}
    The combination of \eqref{eq:1} and \eqref{eq:2} proves {\slshape (ii)}.\vspace{5pt}

    \noindent {\slshape (iii)} This is straightforward using Remark~\ref{CylinderLikeSufficientCondition}.
\end{pr}

\begin{prop}
    The following statements hold:
    \begin{itemize}
        \item[(i)] If $f, g: \R^N \rightarrow [0, \infty)$ are cylinder-like convex functions, then so is $\max(f, g)$;
        \item[(ii)] If $M > N$ are positive integers, $g: \R^N \rightarrow [0, \infty)$ is a cylinder-like convex function and $p: \R^M \rightarrow \R^N$ is projection onto any $N$ distinct coordinates of $\R^M$, then $g \circ p$ is a convex cylinder-like function.
    \end{itemize}

\label{BuildingCylinderLikeFunctions}
\end{prop}

\begin{pr} 
{\slshape (i)}
Call \( h = \max(f,g) \) and note that \( h \) is convex and takes values in \([0, \infty)\). Then, take \(\bar{x}, \bar{y} \in \R^N\) such that \( \bar{x} \neq \bar{y} \) and \( h(\bar{x}) = h(\bar{y}) \). By Remark~\ref{CylinderLikeSufficientCondition} we may assume that $(\bar{x}, \bar{y}) \in H_h^c$, hence $h$ is constant on the line segment between $\bar{x}$ and $\bar{y}$.

We divide the problem into three cases.\vspace{5pt}

\noindent \underline{Case 1}: Assume one of the following holds:
\begin{itemize}
    \item $f(\bar{x}) < h(\bar{x})$ and $f(\bar{y}) < h(\bar{y})$,
    \item $g(\bar{x}) < h(\bar{x})$ and $g(\bar{y}) < h(\bar{y})$.
\end{itemize}
Without losing generality we assume to be under the first assumption, the proofs in the other case is analogous. Then, we have $g(\bar{x}) = h(\bar{x}) = h(\bar{y}) = g(\bar{y})$. By continuity of $f$ and $h$, we deduce that there is $\delta > 0$ such that
\begin{equation*}
    h = g \quad \text{    in } B_\delta(\bar{x}) \cup B_\delta(\bar{y})
\end{equation*}
We can then use the cylinder-like property of $g$ to conclude the proof in this case. \vspace{5pt}

\noindent \underline{Case 2}: Assume one of the following holds:
\begin{itemize}
    \item $f(\bar{x}) = h(\bar{x})$ and $f(\bar{y}) < h(\bar{y})$
    \item $g(\bar{x}) = h(\bar{x})$ and $g(\bar{y}) < h(\bar{y})$
    \item $f(\bar{x}) < h(\bar{x})$ and $f(\bar{y}) = h(\bar{y})$
    \item $g(\bar{x}) < h(\bar{x})$ and $g(\bar{y}) = h(\bar{y})$
\end{itemize}
Suppose that the first assumption holds, the proofs for the other cases are similar. In particular, we have $f(\bar{x}) = h(\bar{x}) = h(\bar{y}) > f(\bar{y})$. By continuity of $f$ along the line segment connecting $\bar{x}$ and $\bar{y}$, we can find $\epsilon_1 > 0$ such that
\begin{equation*}
    f(\bar{x} - \epsilon_1 (\bar{x} - \bar{y})) < f(\bar{x}),
\end{equation*}
Then, the continuity of the function
\begin{equation*}
    z \mapsto f(z - \epsilon_1 (\bar{x} - \bar{y})) - f(z)
\end{equation*}
and of $h$, gives us $\delta_1 > 0$ such that
\begin{equation}
    f(x - \epsilon_1 (\bar{x} - \bar{y})) < f(x) \leq h(x) \quad \forall x \in B_{\delta_1}(\bar{x}).
    \label{maxOfCL: f at x}
\end{equation}
Moreover, as $f(\bar{y}) < f(\bar{x})$, we can find $\rho > 0$ such that
\begin{equation}
    h(y) = g(y) \quad \forall y \in B_{\rho}(\bar{y}).
    \label{maxOfCL: f at y}
\end{equation}

On the other hand, by monotonicity of difference quotients of $f$ along the line segment connecting $\bar{x}$ and $\bar{y}$, for all $t \in (0, 1)$ we have
\begin{align*}
    f((1-t)\bar{x} + t\bar{y}) &< f(\bar{x}) = h(\bar{x}) = h((1-t)\bar{x} + t\bar{y}).
\end{align*}
So we must have $g = h$ on the open segment between $\bar{x}$ and $\bar{y}$. By continuity the equality holds on the closed segment, hence $g(\bar{x}) = g(\bar{y})$. As $g$ is cylinder-like, we can find $\delta_2 > 0$ and $\epsilon_2 > 0$ such that
\begin{equation}
\begin{split}
    g(x - \epsilon_2 (\bar{x} - \bar{y})) \leq g(x) \leq h(x) \quad \forall x \in B_{\delta_2}(\bar{x}),\\
    g(y + \epsilon_2 (\bar{x} - \bar{y})) \leq g(y) \leq h(y) \quad
    \forall y \in B_{\delta_2}(\bar{y}).
\end{split}
\label{maxOfCL: g at both ends}
\end{equation}
Now, defining $\delta := \min (\delta_1, \delta_2)$ and $\epsilon := \min (\epsilon_1, \epsilon_2)$ and keeping in mind Remark~\ref{CylinderLikeSmallerEpsilonRemark}, the combination of  \eqref{maxOfCL: f at x} and \eqref{maxOfCL: g at both ends} guarantees that 
\begin{equation*}
    h(x - \epsilon (\bar{x} - \bar{y})) \leq h(x) \quad \forall x \in B_{\delta}(\bar{x}).
\end{equation*}
Moreover, up to taking smaller $\delta$ and $\varepsilon$ (cf. Remark~\ref{CylinderLikeSmallerEpsilonRemark}), we can ensure that 
\begin{equation*}
    y + \epsilon (\bar{x} - \bar{y}) \in B_{\rho}(\bar{y}) \quad \forall y \in B_{\delta}(\bar{y})
\end{equation*}
and thus the combination of \eqref{maxOfCL: f at y} and \eqref{maxOfCL: g at both ends} yields that
\begin{equation*}
    h(y + \epsilon (\bar{x} - \bar{y})) \leq h(y) \quad \forall y \in B_{\delta}(\bar{y}),
\end{equation*}
concluding the proof of this case.\vspace{5pt}

\noindent \underline{Case 3}: Assume that $f(\bar{x}) = h(\bar{x}) = g(\bar{x})$ and $f(\bar{y}) = h(\bar{y}) = g(\bar{y})$. Let $\delta_1$ and $\epsilon_1$ be the parameters with which \eqref{eq:CL} is satisfied for $f$ and $\delta_2$ and $\epsilon_2$ the ones for $g$. Then, taking $\delta = \min(\delta_1, \delta_2)$ and $\epsilon = \min(\epsilon_1, \epsilon_2)$ and keeping in mind Remark~\ref{CylinderLikeSmallerEpsilonRemark}, it is straightforward to see that \eqref{eq:CL} holds for $h$ with $\delta$ and $\epsilon$.\vspace{5pt}

{\slshape (ii)} First, observe that $g \circ p$ is clearly convex and non-negative. Then, take any distinct $\bar{x}, \bar{y} \in \R^N$ with $g(p(\bar{x})) = g(p(\bar{y}))$. We divide the proof in two cases.

\noindent \underline{Case 1}: Suppose that $p(\bar{x}) \neq p(\bar{y})$.
Then, as $g$ is cylinder-like, we can find $\epsilon > 0$ and $\delta > 0$ such that for all $x \in B_\delta(p(\bar{x}))$ and all $y \in B_\delta(p(\bar{y}))$ we have
\begin{equation*}
\begin{split}
    g(x - \epsilon(p(\bar{x}) - p(\bar{y}))) \leq g(x) \quad \text{and}\quad g(y + \epsilon(p(\bar{x}) - p(\bar{y}))) \leq g(y).
\end{split}
\end{equation*}
Now, for all $\tilde{x} \in B_\delta(\bar{x})$, we have $p(\tilde{x}) \in B_\delta(p(\bar{x}))$, so we can set $x = p(\tilde{x})$ in the above to get
\begin{equation*}
    g(p(\tilde{x} - \epsilon(\bar{x} - \bar{y}))) \leq g(p(\tilde{x})).
\end{equation*}
Similarly, for all $\tilde{y} \in B_\delta(\bar{y})$, we obtain
\begin{equation*}
    g(p(\tilde{y} + \epsilon(\bar{x} - \bar{y}))) \leq g(p(\tilde{y})).
\end{equation*}

\noindent \underline{Case 2}: Suppose that $p(\bar{x}) = p(\bar{y})$.
In this case $g(p(z \pm \epsilon(\bar{x} - \bar{y}))) = g(p(z))$ for all $z \in \R^N$, so the result is trivial. 

The combination of Cases 1 and 2 proves {\slshape (ii)}.
\end{pr}


\begin{example}
By (2), (3) of Proposition~\ref{CylinderLikeFunctionsProp} and all of Proposition~\ref{BuildingCylinderLikeFunctions}, the ``cylinder norm'' of $\R^3$, defined by
\begin{equation*}
    \norm{(x, y, z)} = \max \{ \sqrt{x^2 + y^2}, |z| \}
\end{equation*}
is cylinder-like.
\end{example}

\begin{example}
    Given a finite set of vectors $V \subseteq \R^N$ such that $\text{span}(V) = \R^N$, the associated crystalline norm is defined by
    \begin{equation*}
        \norm{x} := \max_{v \in V} |\scal{x}{v}|.
    \end{equation*}
    Combining {\slshape (i)} of Proposition~\ref{CylinderLikeFunctionsProp} and {\slshape (i)} of Proposition~\ref{BuildingCylinderLikeFunctions}, we deduce that every crystalline norm is cylinder-like.
\end{example}

\noindent We remark that, according to point {\slshape (ii)} of Proposition \ref{CylinderLikeFunctionsProp}, all norms of $\R^2$ are cylinder-like. Unfortunately, this does not hold for higher dimension, as shown by the following example.

\begin{example}
    Define the ``double cone'' norm of $\R^3$ by
    \begin{equation*}
        \norm{(x, y, z)} := \sqrt{x^2 + y^2} + |z|.
    \end{equation*}
    To see this is not cylinder-like, take $\bar{x} = (2/3, 0, 1/3)$ and $\bar{y} = (1/3, 0, 2/3)$ in Definition~\ref{CylinderLike}.
\end{example}

\section{The secondary variational problem}\label{sec:secvar}

In the next section we will prove the existence of optimal transport maps for cost functions of the form $c(x,y) = f(g(y - x))$, where $f: [0, \infty) \rightarrow [0, \infty)$ is increasing and strictly convex, and $g: \R^N \rightarrow [0, \infty)$ is a convex cylinder-like function. As the resulting cost function $c$ is not strictly convex, we cannot expect to have uniqueness of optimal transport maps. However, we are able to prove their existence, solving the secondary variational problem we introduce in this section.

Consider a continuous cost function $c: \R^N \times \R^N \rightarrow [0, \infty)$. Given two measures $\mu, \nu \in \Prob(\R^N)$ such that $\mathscr C(\mu, \nu) < \infty$, consider the usual Kantorivich problem
\begin{equation*}
    \min_{\pi \in \mathsf{Adm}(\mu, \nu)} \int_{\R^N \times \R^N} c(x, y) \de \pi(x, y),
\end{equation*}
and call $\Pi_1(\mu, \nu)$ the set of minimizers. Then, consider the secondary variational problem
\begin{equation}\label{eq:2ndvarproblem}
    \min_{\pi \in \Pi_1(\mu, \nu)} \int_{\R^N \times \R^N} d_{eu}^2(x,y) \de \pi(x, y),
\end{equation}
where $d_{eu}$ denotes the Euclidean distance. Call $\Pi_2(\mu, \nu) \subseteq \Pi_1(\mu, \nu)$ the set of minimisers, which by weak compactness of $\Pi_1(\mu, \nu)$ is not empty. The next statement provides a second order cyclical monotonicity property of the plans in $\Pi_2(\mu, \nu)$. We refer the reader to \cite[Proposition~3.1]{magnabosco3} for the proof.

\begin{prop}\label{prop:doublemonotonicity}
If $\mu$ and $\nu$ have compact support, every $\pi\in \Pi_2(\mu,\nu)$ is concentrated in a set $\Gamma$, such that for every $(x,y),(x',y') \in \Gamma$ it holds that
\begin{equation}\label{eq:monot1}
    c(x,y) + c(x',y') \leq c(x,y') + c(x',y),
\end{equation}
moreover, if $c(x,y) + c(x',y') = c(x,y') + c(x',y)$, then 
\begin{equation}\label{eq:monot2}
    d^2_{eu}(x,y)+d^2_{eu}(x',y') \leq d^2_{eu} (x,y')+ d^2_{eu}(x',y).
\end{equation}
\end{prop}

\section{Proof of the main result}

We can now prove one of the main results of this work.

\begin{theorem} \label{theorem: CylinderLike}
    Let $f: [0, \infty) \rightarrow [0, \infty)$ be an increasing strictly convex function, $g: \R^N \rightarrow [0, \infty)$ be a cylinder-like convex function and consider the cost function $c(x, y) = f(g(y - x))$. Let $\mu, \nu \in \Prob(\R^N)$ with $\mu \ll \Leb^N$ such that $\mathscr C(\mu, \nu) < \infty$.
    Then $\Pi_2(\mu, \nu)$ has a unique element and it is induced by a map.
\end{theorem}

\begin{pr}
    Reasoning as in the proof of Corollary \ref{cor:strictconv}, it suffices to show that every element of $\Pi_2(\mu, \nu)$ is induced by a map. Moreover, according to \cite[Lemma~2.9]{LMT}, we can assume without loss of generality that both $\mu$ and $\nu$ have compact support. Thus take any $\pi \in \Pi_2(\mu, \nu)$. Applying Proposition~\ref{prop:doublemonotonicity}, we can find a full $\pi$-measure set $\Gamma$, satisfying the cyclical monotonicity requirements \eqref{eq:monot1} and \eqref{eq:monot2}.

    Assume by contradiction that $\pi$ is not induced by a map. Let $(\pi_z)_{z \in \R^N}$ be the disintegration kernel of $\pi$ with respect to $\p_1$. Then the set 
    \begin{equation*}
        E := \{ z \in \R^N : \pi_z \text{  is not a delta measure} \}
    \end{equation*}
    has positive $\mu$-measure.
    For any fixed $z \in \R^N$ define the probability measure $\tilde{\pi}_z$ by $\tilde{\pi}_z(B) = \pi_z(B + z)$ for all $B \subseteq \R^N$ Borel.
    Introduce the following probability measure on $\R^N \times \R^N$
     \begin{equation*}
         \gamma := \int_{\R^N} \tilde{\pi}_z \times \tilde{\pi_z} \de \mu(z).
     \end{equation*}
Observe that $\gamma$ is not concentrated on the diagonal set $\{(x, y) \in \R^N \times \R^N : x = y\}$. Therefore, we can take $(\bar{x}, \bar{y}) \in \supp \gamma$, with $\bar{x} \neq \bar{y}$.

     Call $h = f \circ g$ and let $L := \{(x, y) \in \R^N \times \R^N : g(x) \neq g(y)\}$. Notice that, since $f$ is strictly convex, we have $L \subseteq H_h$. Thus, applying Theorem~\ref{theorem: StrictlyConvexCase} to the convex function $h$, we obtain 
     \begin{equation} \label{eq: noMassSplit1}
         \mu(\{z \in \R^N : \supp \tilde{\pi}_z \text{  contains  } x, y \text{  with  } (x, y) \in L\}) = 0.
     \end{equation}
     Since $L$ is open, we have that $(\bar{x}, \bar{y}) \notin L$. Indeed, if this was not the case, we could find $s > 0$ such that
     \begin{equation*}
         \gamma(B_s(\bar{x}) \times B_s(\bar{y})) > 0 \text{  and  } B_s(\bar{x}) \times B_s(\bar{y}) \subseteq L
     \end{equation*} 
     and so
     \begin{equation*}
         \mu(\{z \in \R^N : \tilde{\pi_z}(B_s(\bar{x}))>0, \tilde{\pi_z}(B_s(\bar{y}))>0\}) > 0,
     \end{equation*}
     contradicting \eqref{eq: noMassSplit1}.

     We have then found that $g(\bar{x}) = g(\bar{y})$. As $g$ is cylinder-like, letting $\bar{v} := \bar{x} - \bar{y}$, there exists $\epsilon > 0$ and $\delta_0 > 0$ such that for all $x \in B_{\delta_0}(\bar{x})$ and all $y \in B_{\delta_0}(\bar{y})$
     \begin{equation} \label{cylinderDecrNorm}
         \begin{split}
             g(x - \epsilon\bar{v}) \leq g(x) \quad \text{and} \quad g(y + \epsilon\bar{v}) \leq g(y).
         \end{split}
     \end{equation}
     As $\eta := \scal{\bar{v}}{\bar{x} - \bar{y}} > 0$, we can choose $0 < \delta < \delta_0$ small enough so that for all $x \in B_\delta(\bar{x})$ and all $y \in B_\delta(\bar{y})$ we have
     \begin{equation*}
         \scal{\bar{v}}{x - y} > \eta / 2.
     \end{equation*}
     With this choice of $\delta$, since $(\bar{x}, \bar{y}) \in \supp \gamma$, we have that
     \begin{equation*}
         0 < \gamma(B_\delta(\bar{x}) \times B_\delta(\bar{y})) = \int_{\R^N} \tilde{\pi}_z(B_\delta(\bar{x}))\tilde{\pi}_z(B_\delta(\bar{y})) \de \mu(z),
     \end{equation*}
     and therefore the set
     \begin{equation*}
         F := \{z \in \R^N : \tilde{\pi}_z(B_\delta(\bar{x})) > 0, \tilde{\pi}_z(B_\delta(\bar{y})) > 0\}
     \end{equation*}
     has positive $\mu$-measure. For any $z \in \R^N$ we write $\Gamma_z := \{w \in \R^N : (z, w) \in \Gamma\}$, then $\pi(\Gamma) = 1$ implies that $\tilde{\pi_z}(\Gamma_z) = 1$ $\mu$-almost everywhere. Therefore, the set
     \begin{equation*}
         G := \{z \in \R^N : \tilde{\pi}_z\big(B_\delta(\bar{x}) \cap (\Gamma_z - z)\big) > 0, \tilde{\pi}_z\big(B_\delta(\bar{y}) \cap (\Gamma_z - z)\big) > 0\}
     \end{equation*}
     has positive $\mu$-measure.

     Take a Lebesgue density point $\bar{z}$ of $G$. In a neighbourhood of $\bar{z}$, we can find $z, z + \varepsilon\bar{v}$ for a suitable $0 < \varepsilon < \min\left\{\epsilon, \frac{\eta}{2\scal{\bar{v}}{\bar{v}}}\right\}$. This means there is $x \in B_\delta(\bar{x})$ and $y \in B_\delta(\bar{y})$ such that $(z, z + x), (z + \varepsilon\bar{v}, z + \varepsilon\bar{v} + y) \in \Gamma$.
     Since $0 < \varepsilon < \epsilon$ and $x \in B_\delta(\bar{x}) \subseteq B_{\delta_0}(\bar{x}), y \in B_\delta(\bar{y}) \subseteq B_{\delta_0}(\bar{y})$, it follows from \eqref{cylinderDecrNorm} that
     \begin{equation}\label{eq:twoineq}
         \begin{split}
             g(x - \varepsilon\bar{v}) \leq g(x) \quad \text{and} \quad g(y + \varepsilon\bar{v}) \leq g(y).
         \end{split}
     \end{equation}
     On the other hand, the assumptions on $f$ imply that $f$ is strictly increasing. So if one of the two inequalities in \ref{eq:twoineq} is strict, then
     \begin{align*}
         c(z + \varepsilon\bar{v}, z + x) + c(z, z + \varepsilon\bar{v} + y)&= f(g(x - \varepsilon\bar{v})) + f(g(y + \varepsilon\bar{v}))\\
         &< f(g(x)) + f(g(y))\\
         &= c(z, z+x) + c(z+\varepsilon\bar{v}, z + \varepsilon\bar{v} + y),
     \end{align*}
     contradicting the condition \eqref{eq:monot1} of Proposition~\ref{prop:doublemonotonicity}. Otherwise, an analogous calculation gives that
     \begin{equation*}
         c(z + \varepsilon\bar{v}, z + x) + c(z, z + \varepsilon\bar{v} + y) = c(z, z+x) + c(z+\varepsilon\bar{v}, z + \varepsilon\bar{v} + y)
     \end{equation*}
     but at the same time we have that
     \begin{align*}
         d_{eu}^2(z + \varepsilon\bar{v}, z + x) + d_{eu}^2(z, z + \varepsilon\bar{v} + y) &=\scal{x - \varepsilon\bar{v}}{x - \varepsilon\bar{v}} + \scal{y + \varepsilon\bar{v}}{y + \varepsilon\bar{v}} \\
         &=\scal{x}{x} + \scal{y}{y} + 2\varepsilon(\scal{\bar{v}}{y - x} + \varepsilon\scal{\bar{v}}{\bar{v}}) \\
         &<\scal{x}{x} + \scal{y}{y} + 2\varepsilon\left(-\frac{\eta}{2} + \frac{\eta}{2}\right) \\
         &=d_{eu}^2(z, z + x) + d_{eu}^2(z + \varepsilon\bar{v}, z + \varepsilon\bar{v} + y),
     \end{align*}
    contradicting the condition \eqref{eq:monot2} of Proposition~\ref{prop:doublemonotonicity}. This concludes the proof.
\end{pr}

 Specialising Theorem \ref{theorem: CylinderLike} to the $L^p$-optimal transport problem for cylinder-like norm, we obtain the following corollary, cf. Corollary \ref{cor:strictconv}.

\begin{corollary} \label{cor: cylinderLikeMainApplication}
    Let $\norm{\cdot}$ be a cylinder-like norm on $\R^N$ and $1 < p < \infty$. Consider $c(x, y) = \norm{y - x}^p$ and let $\mu, \nu \in \Prob(\R^N)$ be Borel probability measures with finite $p$-th order moment on $\R^N$, such that $\mu \ll \Leb^N$. 
    Then $\mathscr C(\mu, \nu) < \infty$ and $\Pi_2(\mu, \nu)$ has a unique element which is induced by a map.
\end{corollary}

\begin{remark}
    For a convex function $g: \R^N \rightarrow [0, \infty)$, the condition that $g$ is cylinder-like in Theorem~\ref{theorem: CylinderLike} can be weakened to the condition below, only requiring small modifications to the proof:

    \noindent There is a (countable) partition $\R^N = \bigcup_{i=1}^\infty A_i$ such that whenever $\bar{x} \neq \bar{y}$ and $g(\bar{x}) = g(\bar{y})$, there is $\delta > 0$ and $\epsilon > 0$ such that
    \begin{equation*}
        \begin{split}
            &g(x - \epsilon(\bar{x} - \bar{y})) \leq g(x) \text{   }\forall x\in B_\delta(\bar{x}) \cap A_i      \\
            &g(y + \epsilon(\bar{x} - \bar{y})) \leq g(y) \text{       }\forall y\in B_\delta(\bar{y}) \cap A_j
        \end{split}
    \end{equation*}
    where $A_i, A_j$ are the equivalence classes of $\bar{x}, \bar{y}$ respectively.
    An example of a function which is not cylinder-like but satisfies the above condition is the norm in $\R^3$ whose unit ball is the closed convex hull of the points
    \begin{equation*}
        \{(0, 0, 1), (0, 0, -1)\} \cup \bigcup_{i=1}^{\infty} \{(\cos(\pi / 2^i), \sin(\pi / 2^i), 0), (-\cos(\pi / 2^i), -\sin(\pi / 2^i), 0)\}.
    \end{equation*}
\end{remark}

\bibliography{bibliography}

\bibliographystyle{abbrv} 

\end{document}